\documentclass[12pt,a4paper]{article}
\usepackage{amssymb}
\usepackage{latexsym}
\usepackage{amsmath}
\usepackage{amscd}

\usepackage{amsfonts, amsthm}
\usepackage{amsxtra, amsopn}

\usepackage[cp1251]{inputenc}
\usepackage[english]{babel}

\def\R{\mathbb R}
\def\C{\mathbb C}
\def\N{\mathbb N}

      \def\dC{{\mathbb C}}

      \def\dR{{\mathbb R}}

\def\cD{{\mathcal D}}      
   \def\cH{{\mathcal H}}

\DeclareMathOperator{\supp}{supp}
\DeclareMathOperator{\im}{Im}
\DeclareMathOperator{\Real}{Re}

\DeclareMathOperator{\sgn}{sgn}
\DeclareMathOperator{\dom}{dom}
\DeclareMathOperator{\ran}{ran}

\def\la{\lambda}
\def\ra{\rightarrow}

\def\Ainf{\mathcal{A}_1}
\def\Ab{\mathcal{A}_0}
\def\Ainfinf{\mathcal{A}_\infty}
\def\Azero{\mathcal{A}_2}

\def\Abbb{\mathcal{A}_b}
\def\loc{{\text{\rm loc}}}

\def\ep{\varepsilon}

\def\slr{\mathcal{S}}

\newcommand{\Skindef}{[\raisebox{0.5 ex}{.},\raisebox{0.5 ex}{.}]}

 \newcommand{\wh}{\widehat}
 \newcommand{\ov}{\overline}
 \newcommand{\lln}{\lambda_{0}}

\theoremstyle{plain}
\newtheorem{theorem}{Theorem}[section]
\newtheorem{corollary}[theorem]{Corollary}
\newtheorem{lemma}[theorem]{Lemma}

\newtheorem{definition}[theorem]{Definition}
\newtheorem{remark}[theorem]{Remark}
\newtheorem{example}[theorem]{Example}
\newtheorem{proposition}[theorem]{Proposition}

\numberwithin{equation}{section}

\newtheorem{assumption}[theorem]{Assumption}

\begin{document}

\title{Spectral properties of singular Sturm-Liouville operators
with indefinite weight $\sgn x$}
%\shorttitle{Singular Sturm-Liouville operators with indefinite weight $\sgn x$}
\author{Illya Karabash and Carsten Trunk}

\date{}

\maketitle

%\keywords{}
%\subclass{}

\begin{abstract}
We consider a singular Sturm-Liouville expression
with the indefinite weight $\sgn x$.
To this expression there is naturally a self-adjoint operator
in some Krein space associated.
We characterize  the local definitizability of this operator
in a neighbourhood of $\infty$. Moreover, in this situation,
the point $\infty$ is
a regular critical point. We construct an operator $A=(\sgn x)(-d^2/dx^2+q)$
with non-real spectrum accumulating
to a real point.
The obtained results are applied to several classes of Sturm-Liouville operators.
\end{abstract}

\section{ Introduction}\label{intro}

We consider the singular Sturm-Liouville differential expression
\begin{equation}\label{e a}
a (y) (x) = (\sgn x) (-y''(x) + q(x) y(x)) ,
\qquad x \in \R ,
\end{equation}
with the signum function as indefinite weight and a real potential $q \in L^1_{loc} (\R)$.
We  assume that \eqref{e a} is in the
limit point case at both $-\infty$ and $+\infty$.
This differential expression is naturally connected with a self-adjoint operator $A$
in the Krein space
$(L^2 (\R), \Skindef)$
(see e.g. \cite{CL89}), where the indefinite inner product $[\cdot,\cdot]$ is defined by
\[
[f,g] = \int_R f\overline{g}\ \sgn x
\ dx , \qquad f,g \in L^2 (\R).
\]

The operator $J:f(x) \mapsto (\sgn x) f(x)$ is a fundamental symmetry in the Krein space
$(L^2 (\R), \Skindef)$. Let us define the operator $L:=JA$. Then $L = -d^2/dx^2 +q$
is a self-adjoint Sturm-Liouville operator in the Hilbert space $L^2 (\R)$.
%Clearly,
%\begin{equation}\label{e L}
%L y = -y''+qy , \qquad y\in \dom (L) (:= %\dom (A)).
%\end{equation}
It was shown in \cite{CL89} that
if $L$ is a non-negative operator in the Hilbert space sense then
$A$ is a definitizable operator with $\infty$ as
a regular critical point.

%In general, the operator $A$ may be not definitizable
%(in Section \ref{s Def-ty} we give a criterion).
%However, the essential spectrum of $A$ is still real and it turns out
%that $A$ is locally definitizable over an appropriate subset of $\C$.
%As a  main result we show the equivalence of the semi-boundedness from below
%of the operator $L$ and the local definitizability of the operator $A$ in a
%neighbourhood of $\infty$.
%A similar result was shown in \cite{B07} for a more general class
%of Sturm-Liouville differential expression but with an extra assumption
%on the spectra of certain associated self-adjoint operators.

In general, the operator $A$ may be not definitizable (in Section
\ref{s Def-ty} we give a criterion). However, under certain
assumptions, $A$ is still locally definitizable over an appropriate
subset of $\C$. It seems that the first result of such type was
obtained in \cite{B07} for the operator $y \mapsto \frac 1w
[(py')' + qy]$ with $w$ as indefinite weight function. Note
that in \cite{B07} $w$ may have many turning points, but rather strong
assumptions on the spectra of certain associated self-adjoint
operators are supposed.

As a  main result we show the equivalence of the semi-boundedness
from below of the operator $L$ and the local definitizability of
the operator $A$ in a neighbourhood of $\infty$.
 Moreover, we give a precise description of the domain of definitizability of
$A$. If $L$ is  semi-bounded from below, we show the existence of a decomposition
$A = \Ainfinf \dot +\Abbb$
such that the operator $\Ainfinf$ is similar to a self-adjoint operator in the
Hilbert space sense and $\Abbb$ is a bounded operator, that is,
the point $\infty$ is
a regular critical point.
Hence, the non-real spectrum of $A$ remains bounded.
But, in contrast to the case of a non-negative operator $L$, now
the non-real spectrum
may accumulate to the real axis. We prove in Section \ref{acc}
the existence
of an even continuous  potential $q$ with a sequence
 of non-real eigenvalues of $A$ accumulating to a real point.
This potential $q$ can be chosen in such a way that $A$
is definitizable over $\overline{\C} \setminus \{0\}$.

Finally, in Section \ref{s SomeCl}, we discuss the spectrum
and the sets of definitizability of $A$ for various classes of potentials
$q$.

Differential operators with indefinite weights appears in many
areas of  physics and applied mathematics
 (see \cite{B85,GvdMP87,KapLekH82,Pyat02} and references therein).
Under certain assumptions such operators are definitizable;
this case was studied extensively
(see \cite{BT07,CL89,CN96,CN98,DahL77,FSh00,FadSh02,F96,KarKos06,KWZM03,Kos06,Par03,Vol96,Z05}
and references therein). In
\cite{B07,B06,BPhT06,KarMFAT00,KarKr05,KarM04,KarM06} certain classes
of differential
 operators that contain definitizable as well as not definitizable
operators were considered.

\textbf{Notation:}
Let $T$ be a
linear operator in a Hilbert space
$\mathfrak{H}$. In what follows $\dom
(T)$, $\ker (T)$, $\ran (T)$ are the
domain, kernel, range of $T$,
respectively. We denote the resolvent
set by $\rho (T)$;
$\sigma(T):=\C\setminus\rho(T)$
stands for the spectrum of $T$. By
$\sigma_p (T)$ the set of eigenvalues
of $T$ is indicated. The discrete spectrum $\sigma_{disc} (T)$ is the set of isolated eigenvalues of
finite algebraic multiplicity; the essential spectrum is  $\sigma_{ess} (T):=
\sigma (T) \setminus \sigma_{disc} (T)$. We denote the indicator function of a set
$S$ by $\chi_{S} ( \cdot )$.

 \section{Sturm-Liouville operators with the indefinite weight $\sgn x$\label{s OpDef}}

\subsection{Differential operators \label{ss Op}} \label{definitionen}

We consider the differential expression
\begin{gather} \label{e ell}
\ell (y) (x) =  -y''(x) +q(x)y(x), \quad x \in \R
\end{gather}
with a real potential $q \in L^1_{loc} (\R)$. Throughout this paper
it is assumed that we have
limit point case at both $-\infty$ and $+\infty$. We set
$$
a (y) (x) = (\sgn x) \left( -y''(x) +q(x)y(x) \right), \quad x \in \R.
$$
Let $\mathfrak{D}$ be the set of all $f \in L^2 (\R)$
such that $f$ and $f^{\prime}$ are absolutely continuous with
$\ell (f) \in L^2 (\R)$.
On $\mathfrak{D}$ we define the operators $A$ and $L$ as follows:
\begin{gather*}
\dom(A) = \dom(L) = \mathfrak{D} ,
\qquad A y = a (y) , \qquad L y =  \ell(y) .
\end{gather*}
We equip $L^2 (\R)$ with the indefinite inner product
\begin{equation} \label{SKindef}
[f,g] := \int_{\R} (\sgn x)f(x) \overline{g(x)} dx, \quad f,g \in L^2 (\R).
\end{equation}
Then $(L^2 (\R), \Skindef)$ is a Krein space (for the definition
of a Krein space and basic notions therein we refer to \cite{AJ89}).
A fundamental symmetry $J$ in $(L^2 (\R), \Skindef)$ is given by
\begin{gather*} %\label{e defJ}
(Jf)(x) = (\sgn x) f (x) , \quad f \in L^2 (\R).
\end{gather*}
Obviously,
$$A=JL$$
holds.

Since the differential expressions $a(\cdot)$ and $\ell (\cdot)$
are in the limit point case both at $+\infty$ and $-\infty$, the operator $L$ is self-adjoint in the Hilbert
space $L^2 (\R)$. As  $A=JL$, the operator $A$ is
self-adjoint in the Krein space $L^2 (\R, \Skindef)$.

\begin{definition} \label{AAA}
We shall say that $A$ is \emph{the operator associated with} the
differential expression $a (\cdot)$.
\end{definition}

\subsection{Titchmarsh-Weyl coefficients \label{ss WTc}}

In the following we denote
by  $\C_{\pm}$ the set $ \{ z \in \C : \pm \mbox{Im}\, z >0 \}$.
Let $c_{\lambda} (x)$ and $s_{\lambda} (x)$ denote the fundamental solutions of  the
equation
\begin{equation} \label{e eqL}
-y''(x)+q(x)y(x)=\lambda y(x), \qquad
x\in \R,
\end{equation}
 which satisfy the following conditions
\[
c_{\lambda} (0)= s_{\lambda}'(0)=1;
\qquad c_{\lambda}^{\prime}(0)=
s_{\lambda} (0)=0.
\]
Since the equation \eqref{e eqL} is limit-point at $+\infty$, the
Titchmarsh-Weyl theory (see, for example, \cite{LevSar90}) states that
there exists a unique holomorphic function $m_+ (\lambda)$,
$\lambda \in \C_+ \cup \C_- $, such that the function $s_{\lambda} (\cdot) - m_+ (\lambda) c_\lambda (\cdot) $ belongs to $L^2(\R_+)$. Similarly, the limit point case at
$-\infty$ yields the fact that there exists a unique holomorphic function
$m_- (\lambda)$, $\lambda \in \C_+ \cup \C_- $, such that
$ s_{\lambda} (\cdot) + m_- (\lambda) c_\lambda (\cdot)\in L^2(\R_-)$.
The function $m_+$
($m_-$) is called \emph{the Titchmarsh-Weyl
m-coefficient for} \eqref{e eqL} on
$\R_+$ (on $\R_-$, respectively).

We put
\begin{gather*}
M_\pm (\lambda) := \pm m_\pm (\pm
\lambda) \label{e Mpm} \ .
\end{gather*}

\begin{definition}
The function $M_+ (\cdot)$ ($M_-
(\cdot)$) is said to be \emph{the
Titchmarsh-Weyl coefficient of the
differential expression } $a(\cdot)$
on $\R_+$ (on $\R_-$).
\end{definition}

It is easy to see that for $\la \in
\C_+ \cup \C_-$ the functions
\begin{equation} \label{e def psi}
\psi_\lambda^\pm (x) :=
\begin{cases}
s_{\pm \lambda} (x) - M_\pm (\la )
c_{\pm \la} (x) , & \qquad x \in \R_\pm \\
0, & \qquad x \in \R_\mp
\end{cases}
\end{equation}
belongs to $L^2(\R)$.
Moreover, the following formula (see \cite{LevSar90}) for the norms of $\psi_\lambda^\pm$ in $L^2 (\R)$ holds true
\begin{equation}\label{e |psi|}
\| \psi_\lambda^\pm (x) \|^2 =  \frac {\im M_\pm
(\lambda)}{\im \lambda} , \qquad \lambda \in \C\setminus\R .
\end{equation}

A holomorphic function $G :\C_+ \cup
\C_- \rightarrow \C$ is called \emph{Nevanlinna function} or
\emph{of class (R)}, see e.g.\ \cite{KacKr68}, if
$G (\overline{\lambda})=
\overline{G(\lambda)}$ and
$ \im\lambda \cdot
\im G(\lambda) \geq 0$
for $\lambda \in \C_+ \cup \C_-$.
It follows easily from \eqref{e |psi|} that the functions $M_+$ and $M_-$ (as well as $m_\pm$ ) belong to the class (R).
Moreover, the functions $M_\pm $ have the following  asymptotic behavior
\begin{gather} \label{e asM}
M_\pm (\la ) = \pm \frac{i}{\sqrt{\pm\la}} + O\left( \frac{1}{|\la |}\right) , \quad
(\la \ra \infty , \
0<\delta<\arg \la < \pi - \delta )
\end{gather}
for $\delta \in (0, \frac{\pi}{2})$,  see \cite{Ev72}.
Here and below $\sqrt{z}$ is the branch
of the multifunction on the complex plane
$\C$ with the cut along $\R_+$, singled out by the condition $\sqrt{-1}=i$.

%We assume that $\sqrt{\la} \geq 0$ for
%$\la \in [0, +\infty)$.

\subsection{The non-real spectrum of $A$} \label{nonreal}

In the following we identify functions $f\in L^2(\dR)$ with elements
$\bigl
(\begin{smallmatrix}f_+\\f_-\end{smallmatrix}\bigr)$, where $f_\pm:=f\upharpoonright_{\dR_\pm}\in L^2(\dR_\pm)$.
Similarly we write $q_\pm :=q\upharpoonright_{\dR_\pm}\in L^1_{\loc}(\dR_\pm)$.
Note that the differential expressions
\begin{equation*}
-\frac{d^2}{dx^2}+q_+\qquad\text{and}\qquad  \frac{d^2}{dx^2}- q_-
\end{equation*}
in $L^2(\dR_+)$ and $L^2(\dR_-)$
are both regular at the endpoint $0$ and in the limit point case at the singular endpoint $+\infty$ and $-\infty$,
respectively. Therefore the operators
\begin{equation*}\label{apm}
A^+_{\min}f_+=-f^{\prime\prime}_++q_+f_+\quad\text{and}\quad A^-_{\min}f_-=f^{\prime\prime}_--q_-f_-
\end{equation*}
defined on
\begin{equation*}\label{domapm}
\dom A^\pm_{\min}=\bigl\{f_\pm\in\cD^{\pm}_{\max}: f_\pm(0)=f_\pm^{\prime}(0)=0\bigr\},
\end{equation*}
with
\begin{equation*}
\begin{split}
\cD_{\max}^+&=\bigl\{f_+\in L^2(\dR_+):f_+,f^\prime_+ \mbox{ absolutely continuous, }
 -f_+^{\prime\prime} +q_+ f_+\in L^2(\dR_+)\bigr\},\\
\cD_{\max}^-&=\bigl\{f_-\in L^2(\dR_-):f_-,f^\prime_-\mbox{ absolutely continuous, }
f_-^{\prime\prime} -q_- f_-\in L^2(\dR_-)\bigr\},
\end{split}
\end{equation*}
are closed symmetric operators in the Hilbert spaces
 $L^2(\dR_+)$ and $L^2(\dR_-)$, respectively, cf.\ \cite{W87,W03},
with deficiency indices $(1,1)$. The adjoint operators $(A_{\min}^{\pm})^*$ in
the Hilbert space
 $L^2(\dR_\pm)$ are the usual maximal operators defined on $\cD_{\max}^\pm$.

We introduce the operators
\begin{equation*}
A^+_{0}f_+=-f^{\prime\prime}_++q_+f_+\quad\text{and}\quad A^-_{0}f_-=f^{\prime\prime}_--q_-f_-
\end{equation*}
defined on
\begin{equation*}
\dom A^\pm_{0}=\bigl\{f_\pm\in\cD^{\pm}_{\max}: f^{\prime}_\pm(0)=0\bigr\},
\end{equation*}
Evidently, $A^\pm_{0}$ are self-adjoint extensions of $A_{\min}^{\pm}$
in the Hilbert spaces
 $L^2(\dR_+)$ and $L^2(\dR_-)$, respectively, cf.\ \cite{W87,W03}.
In the following we consider $\dom A^\pm_{\min}$ as subsets of
$L^2(\dR)$. Then above considerations imply the following lemma.

\begin{lemma}
Let $\dom A_{\min}:=
\dom A^+_{\min} \oplus \dom A^-_{\min}$ and let the operator $A_{\min}$ be defined on $\dom A_{\min}$,
$$
A_{\min} := \left(\begin{array}{cc}
A^+_{\min}  & 0\\
0 & A^-_{\min}
\end{array}
\right),
$$
with respect to the decomposition $L^2(\R) = L^2(\R_+) \oplus L^2(\R_-)$.
Then $A_{\min}$ is a closed symmetric operator in the Hilbert space $L^2(\R)$
with deficiency indices $(2,2)$. Moreover, we have
$$
A_{\min} = A\upharpoonright_{\dom A_{\min}}, \qquad
A = A_{\min}^*\upharpoonright_{\mathfrak{D}},
$$
where
\begin{equation*}\label{e domA}
\begin{split}
&\mathfrak{D} = \dom (A) = \\
&= \left\{ f =
\bigl
(\begin{smallmatrix}f_+\\f_-\end{smallmatrix}\bigr)
\in \dom (A_{\min}^+)^*  \oplus \dom (A_{\min}^-)^* :
f_+(0) = f_-(0), f^{\prime}_+(0) = f^{\prime}_-(0) \right\}.
\end{split}
\end{equation*}
\end{lemma}

In the following proposition we collect some spectral properties
of $A$.
\begin{proposition}\label{resolv}
Let $A$ be the operator associated with the differential expression $a(\cdot)$. Then:
\begin{description}
\item[(i)] $ \{ \la \in \C \setminus \R \ : \ M_+ (\la) = M_- (\la ) \} = \sigma_p (A) \setminus \R; $
\item[(ii)] $ \{ \la \in \C \setminus \R \ : \ M_+ (\la) \neq M_- (\la ) \} =  \rho (A) \setminus \R;$
\item[(iii)] $\rho (A) \neq \emptyset $.
\item[(iv)]  The essential spectrum $ \sigma_{ess} (A)$ of $A$ is real and
$$\sigma_{ess} (A) = \sigma_{ess} (A_0^+) \cup \sigma_{ess} (A_0^-). $$
The sets  $\sigma_p(A)\cap \C_{\pm}$ are at most countable with possible
limit points belonging to $\sigma_{ess} (A) \cup \{ \infty \}$.
%Moreover, $\la_0 \in \sigma_p (A)\cap \C_{\pm}$ if and only if $M_+ (\la_0) = M_- (\la_0 ).$
%In   the latter
%case $\dim \mathfrak{L}_{\lambda_0}(A)=\m(\la_0),$ where $\m(\la_0)$ is the multiplicity of
%$\la_0 $ as a zero of the analytic function  $M_+ (\la ) - M_- (\la );$
\end{description}
\end{proposition}
For a proof of Proposition~\ref{resolv} we refer to
\cite[Proposition 2.5]{KarM06} and \cite{KarDis,KarKr05}. We mention only that
the statements \textbf{(iii)} and  \textbf{(iv)}  follow from the first and second statement and~(\ref{e asM}).

\section{Criterions for definitizability \label{s Def-ty}}

\subsection{Definitizable and locally definitizable operators} \label{sectionDEFF}

Let $({\mathcal H},\Skindef)$ be a Krein space and let
$A$ be a closed operator in ${\mathcal H}$. We define the extended
spectrum $\sigma_{e}(A)$ of $A$ by $\sigma_{e}(A):= \sigma(A)$
if $A$ is bounded and
$\sigma_{e}(A):=\sigma(A) \cup \{\infty \}$ if $A$ is unbounded.
We set $\rho_{e}(A) := \ov{\C} \setminus \sigma_{e}(A)$.
A point
$\lambda_{0} \in \C$ is said to belong to the {\em approximative point
spectrum\/} $\sigma_{ap}(A)$ of $A$ if there exists a
sequence $(x_{n}) \subset \dom(A)$  with $\|x_{n} \| =1$,
$n= 1,2, \ldots,$ and $\|(A - \lln)x_{n} \| \to 0$ if
$n \to \infty$. For a self-adjoint operator $A$
in ${\mathcal H}$ all real spectral points of $A$ belong to
$\sigma_{ap}(A)$ (see e.g.\ \cite[Corollary VI.6.2]{B}).

First we recall the notions of spectral points of positive and negative type.

The following definition was given in \cite{LMM1},
\cite{LMaM} (for bounded
self-adjoint  operators).

\begin{definition}  \label{definition++}
For a  self-adjoint operator $A$ in
${\mathcal H}$
a point $\lln \in \sigma(A)$ is called a
spectral point of {\em positive\/} {\rm (}{\em negative\/}{\rm )}
 {\em type of $A$\/}
if $\lln \in \sigma_{ap}(A)$ and for every sequence
$(x_{n}) \subset \dom(A)$ with
$\|x_{n} \| =1$ and $\| (A -\lln)x_{n} \| \to 0$ for
$n \to \infty$, we have
\begin{displaymath}
\liminf_{n \to \infty}\, [x_{n},x_{n}] >0 \;\;\;\; \mbox{{\rm
 (}resp.\ } \limsup_{n \to \infty}\, [x_{n},x_{n}] <0 \mbox{{\rm )}}.
\end{displaymath}
The point $\infty$ is said to be
 of {\em positive\/} {\rm (}{\em negative\/}{\rm )}
 {\em type of $A$\/} if $A$ is unbounded and
 for every sequence
$(x_{n}) \subset \dom(A)$ with
$\lim_{n \to \infty} \|x_{n} \| =0$ and $\| Ax_{n} \| = 1$  we have
\begin{displaymath}
\liminf_{n \to \infty}\, [Ax_{n},Ax_{n}] >0 \;\;\;\; \mbox{{\rm
 (}resp.\ } \limsup_{n \to \infty}\, [Ax_{n},Ax_{n}] <0 \mbox{{\rm )}}.
\end{displaymath}
We denote the set of all points of $\sigma_{e}(A)$ of positive
{\rm (}negative{\rm )} type by $\sigma_{++}(A)$
 {\rm (}resp.\ $\sigma_{--}(A)${\rm )}.
  We shall say that an open subset
$\delta$ of $\overline{\R}$ $(=\R \cup \infty)$ is of {\em positive type\/}
{\rm (}{\em negative type\/}{\rm )} with respect
to $A$ if
\begin{displaymath}
\delta \cap \sigma_{e}(A) \subset \sigma_{++}(A)\;\;\;\; \mbox{{\rm
 (}resp.\ } \delta \cap \sigma_{e}(A)  \subset \sigma_{--}(A){\rm )}.
\end{displaymath}
An open set $\delta$ of $\overline{\R}$ is called of {\em definite type\/} if
$\delta$ is of positive or negative type with respect to $A$.
\end{definition}

The sets $\sigma_{++}(A)$ and $\sigma_{--}(A)$ are contained in $\ov{\R}$.
The non-real spectrum of $A$ cannot accumulate at a point belonging to an
open set of definite type.

Recall, that a self-adjoint operator $A$ in a Krein space
$({\mathcal H}, \Skindef)$ is called
definitizable if
$\rho(A) \ne \emptyset$ and there exists a
rational function $p \ne 0$ having poles only in $\rho(A)$
such that $[p(A)x,x] \geq 0$ for all $x \in {\mathcal H}$.
Then the non-real part of the spectrum of $A$ consists of
no more than a finite number of points. Moreover, $A$ has
a spectral function $E$ defined
on the ring generated by all connected subsets of $\ov{\R}$
whose endpoints do not coincide with the points of some
finite set which is contained in $\{ t \in \R : p(t) = 0 \}
\cup \{\infty \}$ (see \cite{Lan82}).

A self-adjoint operator in a Krein space  is definitizable if
and only if it is definitizable over $\ov{\C}$ in the sense
of the following definition
(see e.g.\ \cite[Definition 4.4]{J6}), which localizes the notion of
definitizability.

\begin{definition}  \label{Deflocdef}
%Let $\Omega$ be a domain in $\ov{\C}$ which is symmetric with 
%respect to $\R$ such that $\Omega \cap \ov{\R} \ne \emptyset$ 
%and $\Omega \cap \C^{+}$ and $\Omega \cap \C^{-}$ are simply 
%connected.
Let $\Omega$ be a domain in $\overline{\C}$ such that
\begin{gather}
\Omega \quad \text{is symmetric with respect to} \ \R, 
\quad \Omega \cap \overline{\R} \neq \emptyset, \label{e Om1}\\
  \text{and the domains} \quad \Omega \cap
\C^{+}, \ \Omega \cap
\C^{-} \quad \text{are simply connected}. \label{e Om2}
\end{gather}
Let $A$ be a self-adjoint
operator in the Krein space $({\mathcal H}, \Skindef)$ such that
$\sigma(A) \cap (\Omega \setminus \ov{\R})$ consists of isolated
points which are poles of the resolvent of $A$, and no point of
$\Omega \cap \ov{\R}$ is an accumulation point of the non-real spectrum
$\sigma(A) \setminus \R$ of $A$.
The operator $A$ is called {\em definitizable over $\Omega$\/},
if the following holds.
\begin{itemize}
\item[{\rm (i)}]
For every
closed subset $\Delta$ of $\Omega \cap \ov{\R}$ there exist
an open neighbourhood ${\mathcal U}$ of $\Delta$ in $\ov{\C}$
and numbers $m \geq 1$,
$M > 0$ such that
\begin{equation} \label{e Sinf}
 \|(A - \lambda)^{-1}\| \leq M (|\lambda | +1)^{2m-2}
 |\mbox{{\rm Im}}\, \lambda |^{-m}
\end{equation}
for all $\lambda \in {\mathcal U} \setminus \ov{\R}$.
\item[{\rm (ii)}]
Every point $\lambda \in \Omega \cap \ov{\R}$
has an open connected neighbourhood $I_{\lambda}$ in $\ov{\R}$ such
that both components of $I_{\lambda} \setminus \{\lambda \}$
 are of definite type
{\rm (}cf.\ Definition {\rm \ref{definition++})} with
respect to $A$.
\end{itemize}
\end{definition}

A self-adjoint operator definitizable over $\Omega$ where $\Omega$ is as in
Definition \ref{Deflocdef} possesses a local spectral function $E$.
For the construction and the properties of
this  spectral function we refer to  \cite{J6}  (see also \cite{J}).
We mention only that $E(\Delta)$ is defined and is a self-adjoint
projection in $({\mathcal H}, \Skindef)$ for every union
$\Delta$ of a finite number of connected subsets $\Delta_{i}$,
$i=1,\ldots ,n$, of $\Omega \cap \ov{\R}$, $\ov{\Delta_{i}}
\subset \Omega \cap \ov{\R}$, such that the endpoints of
$\Delta_{i}$ belong to intervals of definite type.
A real point
$\lambda \in \sigma(A) \cap \Omega$ belongs to
$\sigma_{++}(A)$ if and only if
there exists a bounded open interval $\Delta \subset \Omega$,
$\lambda \in \Delta$,
such that $E(\Delta){\mathcal H}$ is a Hilbert space
(cf.\ \cite{AJT}).
A point $t \in \overline{\R} \cap \Omega$ is called a {\it critical point}
of $A$ if there is no
open subset $\Delta \subset \Omega$ of definite type  with
$t \in \Delta$. The set of critical
points of $A$ is denoted by $c(A)$.\label{cA}
A critical point $t$ is called {\it regular} if there exists
an open deleted neighbourhood $\delta_{0}\subset \Omega$ of $t$ such that
the set of the projections $E(\delta)$ where $\delta$ runs
through all intervals $\delta$ with $\overline{\delta} \subset
\delta_{0}$ is bounded. The set of regular critical points of
$A$ is denoted by $c_{r}(A)$. The elements of $c_{s}(A):=
c(A) \setminus c_{r}(A)$ are called {\it singular} critical points.

We will make use of the following perturbation result,
see \cite{B06}.

\begin{theorem} \label{t FRCrLocDef}
Let $T_1$ and $T_2$ be self-adjoint operators in the Krein space $\cH$, let
$\rho (T_1) \cap \rho (T_2) \cap \Omega
\neq \emptyset$ and assume that
\[
(T_1 - \la_0 I)^{-1} - (T_2 - \la_0 I)^{-1}
\]
is a finite rank operator
for some $\la_0 \in \rho (T_1) \cap \rho (T_2)$. Then $T_1$ is definitizable over $\Omega$ if and only if $T_2$ is definitizable over $\Omega$.

Moreover, if $T_1$ is definitizable over $\Omega$ and
$\Delta \subset \Omega \cap \overline{\R}$ is an open interval with end point
$\eta \in \Omega\cap \overline{\R}$ and $\Delta$ is of positive type
(negative type) with respect to $T_1$, then there exist open interval $\Delta'$,
$\Delta' \subset \Delta$, with endpoint $\eta$ such that $\Delta'$ is of positive type (resp. negative type) with respect to $T_2$.
\end{theorem}

\subsection{Definitizability of $A$}

In this section we will give conditions which ensures the
definitizability of the operator $A$ from Definition \ref{AAA}.
The following definition is needed below.

\begin{definition}
We shall say that the sets $S_1$ and $S_2$ of real numbers are separated by a finite number of points if
there exists a finite ordered set
$\{ \alpha_j \}_{j=1}^{N} $, $N \in \N$,
\[
-\infty = \alpha_0 < \alpha_1 \leq  \cdots
\leq \alpha_{N} < \alpha_{N+1} = +\infty ,
\]
such that one of the sets $S_j$, $j=1,2$,
is a subset of
$\
\bigcup\limits_{ k \ \text{is even} } [\alpha_{k} , \alpha_{k+1} ] $ and another one is a subset of $\
\bigcup\limits_{ k \ \text{is odd} } [\alpha_{k} , \alpha_{k+1} ] $.
\end{definition}

The operator $A^+_{0} \oplus A^-_{0}$, where $A^{\pm}_0$ are defined as in Section \ref{nonreal},
is fundamentally reducible (cf.\ \cite[Section 3]{J88}) in the Krein space
$(L^2 (\R), \Skindef)$ (cf.\ (\ref{SKindef})). Hence the following lemma is a
easy consequence of Definitions \ref{definition++} and \ref{Deflocdef}.
\begin{lemma} \label{l A+Asi+-}
Let $\lambda \in \R$. Then
$\la \in \sigma_{++} (A_0^+ \oplus A_0^-)$ \ ($\la \in \sigma_{--} (A_0^+ \oplus A_0^-)$) if and only if $\la \in \sigma (A_0^+) \setminus \sigma (A_0^-)$
\ ($\la \in \sigma (A_0^-) \setminus \sigma (A_0^+)$, resp.).
The operator $A^+_{0} \oplus A^-_{0}$ is definitizable if and only if the sets $ \sigma(A_{0}^+)$ and
$ \sigma(A_{0}^-) $ are separated by a finite number of points.
\end{lemma}

It follows from Proposition \ref{resolv} and $\sigma (A_0^+ \oplus A_0^-) \subset \R$ that
 $\rho (A) \cap \rho (A_0^+ \oplus A_0^-) \neq \emptyset$. Let \\ 
$\la_0 \in
\rho (A) \cap \rho (A_0^+ \oplus A_0^-)$.
The operators $A_0^+ \oplus A_0^-$ and $A$ are extensions of $A_{\min}$ and \
\[\dim \left( \dom (A_0^+ \oplus A_0^-) /
\dom (A_{\min}) \right) =
\dim \left( \dom (A) /
\dom (A_{\min}) \right) =2 .\]
This implies that
$ (A_0^+ \oplus A_0^- - \lambda_0 I)^{-1}
- (A - \lambda_0 I)^{-1}
$
is an operator of rank 2.
Then \cite{JLan79} and Lemma \ref{l A+Asi+-} imply the following theorem.

\begin{theorem} [\cite{KarDis,KarKr05}] \label{c CrDefSL}
The operator $A$ is definitizable if and only if the sets $ \sigma(A_{0}^+)$ and
$ \sigma(A_{0}^-) $ are separated by a finite number of points.
\end{theorem}

\begin{example} \label{e const1}
Let $q$ be a constant potential, $q(x) \equiv c$, \ $c \in \R$.
It is easy to calculate that $\sigma (A_0^+) = [c,+\infty)$ and
$\sigma (A_0^-) = (-\infty,-c]$.
Thus, Corollary \ref{c CrDefSL} implies that the operator $(\sgn x) (-d^2/dx^2 + c)$ is definitizable in the Krein space
$L^2 (\R, \sgn x \, dx)$ if and only if $c \geq 0$.
\end{example}

\subsection{Local definitizability of $A$ %in neighbourhoods of $\infty$
}
\label{ss ldinf}

In this subsection we consider Sturm-Liouville operators defined as in Section \ref{s OpDef} and we prove
 that the operator $A$ is a definitizable operator in a certain neighbourhood of $\infty$
(in the sense of the Krein space $(L^2 (\R), \Skindef)$) if and only if the operator $L$ is semi-bounded from below
(in the sense of the Hilbert space $L^2 (\R$)).
\begin{remark}
Clearly, $L \geq \eta_0 > -\infty$ whenever $q(x) \geq \eta_0 > -\infty$, $x \in \R$.
\end{remark}

The operator $A_0^+ \oplus A_0^-$ is a self-adjoint operator both in the Hilbert space $L^2 (\R)$ and
in the Krein space $(L^2(\R), \Skindef )$, cf. (\ref{SKindef}).

\begin{lemma} \label{l A+ALocDef}
The following statements are equivalent:
\begin{description}
\item[(i)] The operator $L$ is semi-bounded from below.
\item[(ii)] There exists $R>0$ such that the operator $A_0^+ \oplus A_0^-$ is definitizable over the domain
$\{ \lambda \in \overline{\C} \ : \ |\lambda | > R \}$.
\end{description}
\end{lemma}

\begin{proof}
$(i) \Rightarrow (ii)$. \ Since $A_0^+ \oplus A_0^-$ is a self-adjoint operator in
the Hilbert space $L^2 (\R)$, we see that 
\begin{gather} \label{e sA0+A0-}
\sigma(A_0^+ \oplus A_0^-) \subset \R \quad \text{and
%inequality 
\eqref{e Sinf} holds for all} \ \la \in \C \setminus \R \quad \text{with} \ m=1.
\end{gather}

Assume that $L \geq \eta_0 $.
The operator $L$ is a self-adjoint extension of $A^+_{\min} \oplus (-A^-_{\min})$, hence
the operator $A_{\min}^+$ is semi-bounded from below,
$A_{\min}^+ \geq \eta_0$,
 and $A_{\min}^-$ is semi-bounded from above, $A_{\min}^- \leq -\eta_0$
The operators $A_0^\pm$ are self-adjoint extensions
in $L^2 (\R_\pm)$ of the symmetric operators $A_{\min}^\pm$
 with deficiency indices (1,1). Hence the spectrum of $A_0^+$
($A_0^-$) lies, with the possible exception of at most one normal eigenvalue,
in $[\eta_0 , \infty)$ (in $(-\infty, -\eta_0]$, respectively),
see e.g. \cite[Section VII.85]{AG}.

Choose $R:= \eta_0$.
Lemma \ref{l A+Asi+-} implies that
the set $(R,+\infty)$, with the possible exception of at most one  eigenvalue, is of positive type
and the set $(-\infty,-R)$,
 with the possible exception of at most one  eigenvalue,
 is of negative type with respect to $A_0^+ \oplus A_0^-$.
Thus, the operator $A_0^+ \oplus A_0^-$ is definitizable over $\{ \lambda \in \overline{\C} \ : \ |\lambda | > R \}$.

$(i) \Leftarrow (ii)$
Obviously, the Sturm-Liouville operator $A_0^+$ ($A_0^-$) is not semi-bounded from above (below, resp.). That is,
\begin{gather} \label{e supA+infA-}
\sup \sigma (A_0^+) = +\infty , \qquad
\inf \sigma (A_0^-) = -\infty .
\end{gather}

Assume that $L$ is not semi-bounded from below. Then  $A_{\min}^+$ or $-A_{\min}^-$ is not semi-bounded from below. Thus,
$ \inf \sigma (A_0^+) = -\infty $
or $ \sup \sigma (A_0^-) = +\infty $.

Consider the case
\begin{gather} \label{e infA+}
\inf \sigma (A_0^+) = -\infty  . \end{gather}
It follows from \eqref{e infA+}, \eqref{e supA+infA-} and Lemma \ref{l A+Asi+-} that \[
(-\infty,-r) \cap \sigma_{++} (A_0^+ \oplus A_0^-) \neq \emptyset \qquad \text{and} \qquad
(-\infty,-r) \cap \sigma_{--} (A_0^+ \oplus A_0^-) \neq \emptyset
\]
for all $r>0$.
Thus, by definition, the operator $A_0^+ \oplus A_0^-$ is not definitizable over
$\{ \lambda \in \overline{\C} \ : \ |\lambda | > r \}$
for arbitrary $r>0$ .
The case $\sup \sigma (A_0^-) = +\infty$ can be considered in the same way.
\end{proof}

The following  theorem is one of the main results.

\begin{theorem} \label{t DefInf}
The following assertions are equivalent:
\begin{description}
\item[(i)] The operator $L$ is semi-bounded from below.
\item[(ii)] There exists $R>0$ such that the
operator $A$ is definitizable over the domain
$\{ \lambda \in \overline{\C} \ : \ |\lambda | > R \}$.
\end{description}
\end{theorem}

\begin{proof}
It follows from Proposition \ref{resolv}~(iii) and $\sigma (A_0^+ \oplus A_0^-) \subset \R$ that $\rho (A) \cap \rho (A_0^+ \oplus A_0^-) \neq \emptyset$. Let $\la_0 \in
\rho (A) \cap \rho (A_0^+ \oplus A_0^-)$.
The operators $A_0^+ \oplus A_0^-$ and $A$ are extensions of $A_{\min}$ and \[
\dim \left( \dom (A_0^+ \oplus A_0^-) /
\dom (A_{\min}) \right) =
\dim \left( \dom (A) /
\dom (A_{\min}) \right) =2 .\]
This implies that
\begin{equation} \label{2007}
(A_0^+ \oplus A_0^- - \lambda_0 I)^{-1}
- (A - \lambda_0 I)^{-1}
\end{equation}
is an operator of rank 2.
Combining Lemma \ref{l A+ALocDef} and
Theorem \ref{t FRCrLocDef}, Theorem~\ref{t DefInf} is proved.
\end{proof}

By Theorem \ref{t DefInf}, the semi-boundedness of $L$ implies the
definitizability of $A$ over some domain. Now %In the next theorem
we give a precise description of the domain of definitizability of
$A$ in terms of the spectra of $A_0^+$ and $A_0^-$.
% we will
%describe the local definitizability of $A$ over a suitable
%neighbourhood of a given point $\lambda \in \R$. We mention that a
%similar assertion holds for $\lambda = \infty$.

Let $T$ be an operator such that $\sigma(T) \subset \R$. Let us
introduce the sets $\sigma^{left} (T)$ and $\sigma^{right} (T)$ by the following way: 
a point $\la \in \overline{\R}$ $(=\R \cup \infty)$ is said to belong to 
$ \sigma^{left} (T)$ \ $(\sigma^{right} (T))$ if
there exists an increasing (resp.\ decreasing) sequence  $\{\la_n \}_1^\infty \subset \sigma(T)$ such that 
$ \lim_{n\to\infty} \la_n = \la$.

Note that 
\begin{gather} \label{e lress}
\sigma^{left} (T) \cup \sigma^{right} (T) \subset \sigma_{ess} (T)\cup \{\infty\} .
\end{gather} 
For differential operators $A_0^\pm$, equality holds in \eqref{e lress} since every point of $\sigma_{ess} (A_0^\pm)$ is an accumulation point of $\sigma (A_0^\pm)$.
%one can easily obtain 
%\[
%\sigma^{left} (A_0^\pm) \cup \sigma^{right} (A_0^\pm) =
%\sigma_{ess} (A_0^\pm) \cup \{\infty\}.
%\]

We put 
\begin{gather} \label{e slr}
\slr_A := \left( \sigma^{left} (A_0^+) \cap
\sigma^{left} (A_0^-) \right) \cup
\left( \sigma^{right} (A_0^+) \cap
\sigma^{right} (A_0^-) \right).
\end{gather}

\begin{theorem} \label{t DomLD}
Let $\Omega$ be a domain in $\overline{\C}$ such that
\eqref{e Om1}-\eqref{e Om2} are fulfilled.
%\begin{gather}
%\Omega \quad \text{is symmetric with respect to} \ \R, 
%\quad \Omega \cap \overline{\R} \neq \emptyset, \label{e %Om1}\\
%  \text{and} \quad \Omega \cap
%\C^{\pm} \quad \text{are simply connected}. \label{e Om2}
%\end{gather}
Then the operator $A= (\sgn x) (-d^2/dx^2 +q)$ is definitizable
over $ \Omega $ if and only if $\Omega \subset \Omega_A$, where
$\Omega_A:= \overline{\C} \setminus \slr_A $.
\end{theorem}

\begin{proof}
Arguments from the proof of Theorem \ref{t DefInf} show that it is enough to prove the theorem for the operator $A_0^+ \oplus A_0^-$.

Let $\la \in \slr_{A}$ and let $I_\lambda$ be an open connected neighbourhood of $\lambda$. Then \eqref{e slr} and Lemma \ref{l A+Asi+-} imply that one of the components of $I_\lambda \setminus \{ \lambda \}$ is not of definite type. So if 
$A_0^+ \oplus A_0^-$ is definitizable over $\Omega$, then $\la \not \in \Omega$.

Conversely, if $ \slr_{A} \neq \ov{\R}$, then condition (ii) from Definition \ref{Deflocdef} is fulfilled for $\Omega_A = \ov{\C} \setminus \slr_A$. Taking \eqref{e sA0+A0-} into account, we see that $A_0^+ \oplus A_0^-$ is definitizable over $\Omega_A$.
\end{proof}

\begin{remark} \label{r GSLD}
Note that $\Omega_A \cap \overline{\R} = \emptyset$ is equivalent to $\sigma_{ess} (A_0^+) = \sigma_{ess} (A_0^-) = \R$. In the converse case, \eqref{e Om1}-\eqref{e  Om2} are fulfilled for $\Omega_A$ and it is the greatest domain over which the operator $A$ is definitizable.
\end{remark}

The following statement is a simple consequence of Theorem \ref{t DefInf}, Theorem \ref{t DomLD}, and \eqref{e lress}.

\begin{corollary} \label{genau Def}
Assume that $L$ is semi-bounded from below.
Then the operator $A$ is definitizable over the set
$\ov{\C} \setminus (\sigma_{ess} (A_0^+) \cap \sigma_{ess} (A_0^-))$.
\end{corollary}

\subsection{Regularity of the critical point $\infty$}

In the sequel we will use a result which follows  easily from
\cite[Lemma 3.5 (iii)]{CL89} and \cite[Theorem 3.6 (i)]{CL89}.

\begin{proposition}%[\cite{CL89}] 
\label{CL89}
If the operator $\widetilde{L}:=-d^2/dx^2+\widetilde{q}(x)$,
for some real $\widetilde{q}\in L^1_{loc} (\R)$,
 defined on $\mathfrak D$ is nonnegative in the Hilbert space
$L^2 (\R)$, then the operator $\widetilde{A}:=(\sgn
x)\widetilde{L}$ is definitizable and $\infty$ is a regular critical point of
$\widetilde{A}$.
\end{proposition}

The following theorem can be considered as the main result of this note.
\begin{theorem} \label{t RegInf}
Assume that assertions (i), (ii) of Theorem \ref{t DefInf} hold true. Then
there exists a decomposition
%$A = \Ainfinf \dot +\Abbb$
\begin{gather} \label{e Adec}
A = \Ainfinf \dot +\Abbb
\end{gather}
such that the operator $\Ainfinf$ is similar to a self-adjoint operator in the Hilbert space sense and $\Abbb$ is a bounded operator.
\end{theorem}

\begin{remark} The conclusion of Theorem \ref{t RegInf} is equivalent to the
regularity of critical point $\infty$ of  the operator $A$.
\end{remark}

\begin{proof}[Proof of Theorem \ref{t RegInf}]
Assume that  $A$ is an operator definitizable  over $ \{ \lambda \in \overline{\C} \ : \ |\lambda | > R
\}$, $R >0$. By Theorem \ref{t DefInf}, this is equivalent to the fact that $L \geq \eta_0$ for certain $\eta_0 \in \R$.

Denote by $E^A$ the spectral function of $A$. Choose $r>R$ such that
$\sigma(A) \setminus \R \subset \{ \lambda \in \C  \ : \ |\lambda | \leq r \}$ and
$E^A(\overline{\R} \setminus (-r,r))$ is defined. Then $A$ decomposes,
\begin{equation*}
\begin{split}
A = \Ainf \dot + \Ab , \qquad&
\Ainf := A \upharpoonright \dom (A) \cap  (E^A(\overline{\R} \setminus (-r,r))L^2(\R)), \\
& \Ab := A \upharpoonright \dom (A) \cap ((I-E^A(\overline{\R} \setminus (-r,r)))L^2(\R))
\end{split}
\end{equation*}
and the following statements holds (cf.\ \cite[Theorem 2.6]{J88}):
\begin{description}
\item $\Ainf$ \quad is a definitizable operator in the Krein space
$(E^A(\overline{\R} \setminus (-r,r))L^2(\R), \Skindef)$;
\item $\Ab$ \quad is a bounded operator
and  \ $\sigma (\Ab) \subset \{\la : |\la | \leq r \}$.
\end{description}

Let us show that $\infty$ is not a singular critical point of $\Ainf$.

Consider the operator $\Azero$
defined by $\Azero = \Ainf \dot + \, 0$, where the direct sum is considered
with respect to the decomposition
$$
L^2 (\R) = E^A(\overline{\R} \setminus (-r,r))L^2(\R)\dot + (I-E^A(\overline{\R} \setminus (-r,r)))L^2(\R),
$$
and $0$ is the zero operator in the subspace $\ran(I-E^A(\overline{\R} \setminus (-r,r)))$.
 Since $\Ab$ is a bounded operator, we have
\begin{gather*} \label{e domA+0=domA}
\dom (\Azero) = \dom A .
\end{gather*}
It is easy to see that $\Azero$ is a definitizable operator in the Krein space $(L^2 (\R), \Skindef)$.
Moreover, $\infty$ is not a singular critical point of $\Azero$ if and only if $\infty$ is not
a singular critical point of $A$.

Now we prove that $\infty$ is not a singular critical point of $\Azero$.
Let $\eta_1 < \eta_0$. Since $L \geq \eta_0$, we see that $L - \eta_1 I$ is a uniformly positive operator in the Hilbert space $L^2 (\R)$ (i.e.,
$L - \eta_1 I \geq \delta >0$). Therefore $\widetilde A := J (L - \eta_1 I)$,
\begin{gather*}\label{e domtA=domA}
\widetilde A y (x) = (\sgn x) (-y'' (x) + q(x) y(x) - \eta_1 y(x)) , \qquad
\dom(\widetilde A) = \dom (A),
\end{gather*}
is a definitizable nonnegative operator in the Krein space $(L^2 (\R), \Skindef)$.
By Proposition \ref{CL89}, $\infty$ is not a singular critical point of $\widetilde A$. The \'Curgus criterion of the regularity of critical point $\infty$, see \cite[Corollary 3.3]{C85}, implies that $\infty$ is not a singular critical point of the operator $\Azero$. So $\infty$ is not a singular critical point of $\Ainf$.

It follows from $L \geq \eta_0$ and
Lemma \ref{l A+Asi+-} that for sufficiently large $r_1 > 0$ the set $( - \infty, - r_1] $ is of negative type and the set $[ r_1 , +\infty) $ is of positive type with respect to $A_0^+ \oplus A_0^-$.
Combining this with Theorem \ref{t FRCrLocDef}, we obtain that there exists $r_2 \geq r_1$ such that $( - \infty, - r_2] $ is of negative type and the set $[ r_2 , +\infty) $ is of positive type with respect to the operator $A$. Evidently,
we obtain the desired decomposition
\begin{equation*}
\begin{split}
A = \Ainfinf \dot +\Abbb, \qquad&
\Ainfinf := A \upharpoonright \dom (A) \cap  (E^A(\overline{\R} \setminus (-r_2,r_2))L^2(\R)), \\
& \Abbb := A \upharpoonright \dom (A) \cap ((I-E^A(\overline{\R} \setminus (-r_2,r_2)))L^2(\R)),
\end{split}
\end{equation*}
where $ \Abbb$  is a bounded operator and
$\Ainfinf$ is similar to a self-adjoint operator in the Hilbert space sense.
\end{proof}

\section{Accumulation of non-real eigenvalues to a real point}
\label{acc}

By Proposition \ref{resolv} (i), the non-real spectrum $\sigma(A) \setminus \R$ of $A$ consists of eigenvalues.

Let $\slr_A$ be the set defined by \eqref{e slr}.
The following proposition is a consequence of 
Theorems \ref{t DomLD} and \ref{t DefInf}. 

\begin{proposition} \label{p nonreal}
If $\la$ is an accumulation point of $\sigma (A) \setminus \R$,
then $\la \in \slr_A$. In particular, if the operator $L=-d^2/dx^2+q(x)$ is semi-bounded from below, then non-real spectrum of $A$
%$\sigma (A) \setminus \R$ 
is a bounded set.
\end{proposition}

The goal of this subsection is to
show that there exists a potential $q$  continuous in $\R$ 
such that the set of non-real eigenvalues of
the operator $A=(\sgn x)(-d^2/dx^2+q(x))$ \emph{has}
a real accumulation point.

It is well known (e.g. \cite{LevSar90}) that $M_+ $, the Titchmarsh-Weyl
m-coefficient for \eqref{e eqL} (see Subsection \ref{ss WTc}), admits the
following integral representation
\begin{gather*} \label{e intM}
M_+ (\lambda ) = \int_{\R} \frac{d \Sigma_+ (t)}{t-\lambda} , \qquad \lambda \in \C \setminus \R ,
\end{gather*}
where $ \Sigma_+ (\cdot)$ is a nondecreasing scalar functions
such that $ \int_{\R} (1+|t|)^{-1} d \Sigma_+ (t) < \infty $.
The function $\Sigma_+ $ is called \emph{a
spectral function of} the boundary value problem
\begin{gather}\label{e SLNC}
-y''(x)+q_+(x)y(x)=\lambda y(x), \quad y'(0)=0, \qquad x\in[0,+\infty) .
\end{gather}
This means that the self-adjoint operator $A_0^+$ introduced in Subsection \ref{nonreal}
is unitary equivalent to the operator of
multiplication by the independent variable in the Hilbert space
$L^2(\R, d\Sigma_+ (t))$.
This fact obviously implies
\begin{equation} \label{e sA=supp}
\sigma(A_0^+)=\supp(d\Sigma_+),
\end{equation}
where $\supp d\tau$ denotes the
\emph{topological support} of a Borel measure $d\Sigma_+$ on $\R$
(i.e., $\supp d\Sigma_+$ is the smallest closed set $\Omega\subset\R$
such that $d\Sigma_+ (\R \setminus \Omega) = 0$).

\begin{lemma} \label{l iSpEven}
Assume that $q$ is an even potential, $q(x) = q(-x)$, $x \in \R$.
If $\ep >0$, then \ $i\ep \in \sigma_p (A) $ if and only if $\Real M_+ (i\ep) = 0$.
\end{lemma}

\begin{proof}
Since $q$ is even, we get $m_+ (\la) = m_- (\la)$, $\la \in \C \setminus \R$.
So $M_- (i\ep) = - M_+ (-i\ep) $. Since $M_+$ is a Nevanlinna function,
we see that $M_+ (-i \ep) = \overline{M_+ (i \ep)}$. Thus,
\[
M_+ (i\ep) - M_- (i\ep) = M_+ (i\ep) + \overline{M_+ (i \ep)} = 2 \Real M_+ (i \ep).
\]
Proposition \ref{resolv} completes the proof.
\end{proof}

The following lemma follows easily from the Gelfand--Levitan theorem
(see e.g. \cite[Subsection 26.5]{Nai68II}).

\begin{lemma} \label{l InvTh}
Let $ \Sigma (t)$, $t \in \R$, be a nondecreasing function such that
\begin{gather}
\int_{-\infty}^{T_1-0} d \Sigma (t)= 0 \qquad \text{and}  \label{e T1}\\
\int_{-\infty}^{s-0} d \Sigma (t) = \int_{0}^s \frac 1{\pi \sqrt{t}} dt \
\left( = \frac 2\pi \sqrt{s} \right) \qquad \text{for all} \quad s > T_2 . \label{e T2}
\end{gather}
with certain constants $T_1,T_2 \in \R$, $T_1 < T_2$.
Then there exists a potential $q_+$ continuous in $[0,+\infty)$ such that $ \Sigma (t)$
is a spectral function of the boundary value problem
\begin{equation*}\label{e SLq0}
-y''(x)+q_+ (x)y(x)=\lambda y(x), \quad y'(0)=0, \qquad x\in[0,+\infty).
\end{equation*}
\end{lemma}

\begin{lemma} \label{l exS}
There exist a nondecreasing function $\Sigma (t)$, $t \in \R$, with the following
properties:
\begin{description}
\item[(i)] $\Sigma (t) = \Sigma_1 (t) + \Sigma_2 (t)$, where
\begin{equation} \label{e S1}
\Sigma_1 \in AC_{loc} (\R), \qquad \Sigma_1^{'} (t) = \left\{
\begin{array}{rr}
                    0, &   \quad t \in (-\infty,1), \\
                  \frac{1}{\pi\sqrt{t}}, &   t \in (1, + \infty),
\end{array} \right.
\end{equation}
and the measure $d\Sigma_2$ has the form
\begin{multline} \label{e S2}
d\Sigma_2 (t) = \sum_{k=1}^{+\infty} h_k \delta (t-s_k), \qquad \\
h_k >0 , \quad s_k \in (-1,1), \quad k \in \N; \qquad \sum_{k=1}^{+\infty} h_k < \infty ,
\end{multline}
(here $\delta(t)$ is the Dirac delta-function).

\item[(ii)]
%$s_k \in (-1,1)$,  $ \quad k \in \N$.
Conditions \eqref{e T1}-\eqref{e T2} are valid for $\Sigma$ with $T_1=-1$ and $T_2=1$.
%\[
%\int_{-\infty}^s d \Sigma =  \int_{0}^s \frac 1{\pi \sqrt{t}} dt \
%\left( = \frac 2\pi \sqrt{s} \right) \qquad \text{for all} \quad s > 1 .
%\]

\item[(iii)] There exists a sequence $\ep_k>0 $, $k \in \N$, such that $\lim_{k \to \infty} \ep_k = 0$ and
$r (\ep_k) = 0$, $k \in \N$, where the function $r (\ep)$, $\ep >0$, is defined by
\[
r(\ep) :=  \Real \int_\R \frac 1{t-i\ep} d \Sigma (t)=
\int_\R \frac t{t^2+ \ep^2} d \Sigma (t).
\]
\end{description}
\end{lemma}

\begin{proof}
Let $h_k = 2^{-k+1}/\pi$. Then
\begin{equation} \label{e hk}
\sum_{k=1}^{\infty} h_k = 2/\pi .
\end{equation}
Now, if $s_k \in (-1,1)$ for all $k \in \N$, then $\Sigma$ possesses property (ii).
We should only choose $\{ s_k \}_1^\infty \subset (-1,1)$ such that statements (iii) holds true.

Consider for $\ep \geq 0$ the functions
\[
r_0 (\ep) =  \int_1^\infty \frac t{t^2+ \ep^2}  d \Sigma_1 (t)
\]
and
\[
r_n (\ep) := \int_1^\infty \frac t{t^2+ \ep^2} d \Sigma_1 (t) + \sum_{k=1}^n \frac {s_k h_k }{s_k^2+ \ep^2} ,
\qquad n \in \N .
\]

Let $s_k \neq 0$ for all $k \in \N$. Then $r_n$ are well-defined and continuous on $[0,+\infty)$.
Besides, \\ $\lim_{n \to \infty} r_n (\ep) = r (\ep)$ for all $\ep >0$.
It is easy to see that $\lim_{\ep \to \infty} r_n (\ep) = 0$, $n \in \N$.
Since $r_n$ are continuous on $[0,+\infty)$, we see that
\[
\mathrm{SUP}_n := \sup_{\ep \in [0,+\infty)} |r_n (\ep)| < \infty , \qquad n \in \N.
\]

Now we give a procedure to choose $s_k \in (-1,1) \setminus \{ 0 \}$.

Let $s_1$ be an arbitrary number in $(-1,0)$ such that
\begin{gather*}
\left. \frac { s_1 h_1}{s_1^2 +\ep^2} \right|_{\ep = |s_1|} = \frac {1}{\pi} \frac {1}{2 s_1}
< - \mathrm{SUP}_0 -1, \\
\text{in other words}, \qquad - \frac{1}{ 2 \pi (\mathrm{SUP}_0+1)} < s_1 <0 .
\end{gather*}
Then
\begin{gather} \label{e r1<}
r_1 (|s_1|)= r_0 (|s_1|) + \left. \frac {s_1 h_1 }{s_1^2 +\ep^2} \right|_{\ep = |s_1|} < r_0 (|s_1|) - \sup_{\ep \in [0,+\infty)} |r_0 (\ep)| -1 < - 1 .
\end{gather}

Let
\begin{gather} \label{e b1}
 \{s_k \}_2^\infty \in (-b_1,b_1) \setminus \{ 0 \} \quad \text{with certain} \quad b_1 \in (0,|s_1|/2).
\end{gather}
Let us show that we may choose a number $b_1 $ such that \eqref{e b1} implies
\begin{gather} \label{e rs1}
r(|s_1|) < 0 .
\end{gather}
Indeed, \eqref{e r1<} and \eqref{e hk} yield
\begin{gather*}
r(|s_1|) = r_1 (|s_1|) + \left[ \sum_{k=2}^{\infty} \frac {s_k h_k }{s_k^2+ \ep^2}\right]_{\ep = |s_1|} < \\
< -1 + \sum_{k=2}^\infty \frac {h_k |s_k|}{s_k^2+ s_1^2} <
-1 + \frac {b_1}{s_1^2} \sum_{k=2}^\infty h_k < -1 + \frac {2b_1}{\pi s_1^2}
\end{gather*}
and therefore \eqref{e rs1} is valid whenever $0<b_1 < \pi s_1^2 /2$.

Similarly, there exist $s_2 \in (0,b_1) $ such that
\[
\left. \frac { s_2 h_2}{s_2^2 +\ep^2} \right| _{\ep = s_2} = \frac {1}{2\pi} \frac {1}{2 s_2}  >
\mathrm{SUP}_1 +1,
\]
and therefore
\begin{gather*}
r_2 (s_2) > 1 .
\end{gather*}
Further, there exist $b_2 \in (0, s_2/2)$ such that
$ % \begin{gather*}
 \{s_k \}_3^\infty \subset (-b_2,b_2) \setminus \{ 0 \}
$ %\end{gather*}
implies that
$ %\begin{gather*}
r(s_2) > 0 .
$ %\end{gather*}

Continuing this process, we obtain a sequence $\{ s_k \}_1^\infty \subset (-1,1) \setminus \{ 0 \}$ with
the following properties:
\begin{gather}
s_k \in (-1,0) \quad \text{if k is odd}, \qquad
s_k \in (0,1) \quad \text{if k is even} , \notag \\
|s_1| > \frac{|s_1|}{2}>|s_2| > \frac{|s_2|}{2} > |s_3|> ... > |s_k| > \frac{|s_k|}{2} >
|s_{k+1}| > ... \ , \label{e s>s} \\
r(|s_k|) <0 \quad \text{if k is odd}, \qquad
r(|s_k|) >0 \quad \text{if k is even} \label{e rs><0}.
\end{gather}
It is easy to show that $r$ is continuous on $(0,+\infty)$. Combining this with \eqref{e rs><0},
we see that there exists $\ep_k \in (|s_{k-1}|, |s_{k}|) $ such that $ r(\ep_k) = 0$, $k \in \N$.
Besides, \eqref{e s>s} implies $\lim |s_k| = \lim \ep_k = 0$.
\end{proof}

\begin{theorem} \label{t ac-0Ex}
There exist an even  potential $\wh q$  continuous  on $\R$ %, $\wh q (x) = \wh q (-x)$, $x \in \R$,
and a sequence $\{ \ep_k \}_1^{\infty} \subset \R_+ $  such that
\begin{description}
\item[(i)] the operator $\wh A$ defined by the differential
expression
\begin{equation} \label{e de wtq}
(\sgn x)\left(-\frac{d^2}{dx^2}+\wh q(x)\right)
\end{equation}
on the natural domain $\mathfrak{D}$ (see Subsection \ref{ss Op})
is self-adjoint in the Krein space $L^2 (\R, \Skindef)$;
\item[(ii)] $ \{ i \ep_k \}_1^{\infty} \subset \sigma_p (\wh A)$, i.e., $\ i \ep_k$, $k \in \N$, are non-real
eigenvalues of $\wh A$;
\item[(iii)] $\lim_{k \to \infty} \ep_k = 0$;
\item[(iv)] the operator $\wh A$ is definitizable over the domain
$ \overline{\C} \setminus \{ 0 \} $.
\end{description}
\end{theorem}

\begin{proof}
\textbf{(i)} Let $\Sigma$ and $\{ \ep_k \}_1^{\infty}$ be from Lemma \ref{l exS}. Then, by Lemma \ref{l InvTh},
$\Sigma$ is a spectral function of the boundary value problem \eqref{e SLNC} with a certain potential
$\wh q_+$. Let us consider an even  continuous  potential $\wh q(x) = \wh q_+ (|x|)$, $x \in \R$,
and the corresponding operator $\wh A=(\sgn x)\left(-\frac{d^2}{dx^2}+\wh q(x)\right)$
defined as in Subsection \ref{ss Op}.

It is well known that if equation \eqref{e eqL} is in the limit-circle case at $+\infty$ then
$M_+ (\cdot)$ is a meromorphic function on $\C$ and the spectral function $\Sigma_+$ is a step function
with jumps at the poles of $M_+ (\cdot)$ only (see e.g. \cite[Theorem 9.4.1]{CL55}).
As $\Sigma_+(t) = \Sigma(t)$, $t>0$, this condition
does not hold for the function $\Sigma$ since $\Sigma$ satisfies \eqref{e T2}.
Indeed, \eqref{e T2} means that $\Sigma^\prime (t) = \frac 1{\pi \sqrt{t}}$ for $t > T_2 =1$
and therefore $\Sigma$ is not a step function.
So \eqref{e eqL} is limit-point at $+\infty$.

Since the potential $\wh q$ is even, the same is true for $-\infty$.
Thus, $\wh A$ is a self-adjoint operator in the Krein space $L^2 (\R, \Skindef)$,
see Subsection \ref{definitionen}.

\textbf{(ii)} and \textbf{(iii)} follow from Lemma \ref{l iSpEven} and statement (iii) of Lemma \ref{l exS}.

\textbf{(iv)} Let $\wh A_0^\pm$ be the self-adjoint operators in
the Hilbert spaces $L^2 (\R_\pm)$ defined by the differential expression
\eqref{e de wtq} in the same way as in Subsection \ref{nonreal} where $q$ is replaced by $\wh q$.
By \eqref{e sA=supp}, $\sigma (\wh A_0^+) = \{ s_k \}_1^\infty \cup [1,+\infty)$. Since
$\wh q$ is even, one gets $\sigma (\wh A_0^-) = \{ -s_k \}_1^\infty \cup (-\infty,-1]$. It follows from
$ \{ s_k \}_1^\infty \subset (-1,1)$ and $\lim_{k \to \infty} s_k = 0$
that
$$
\min \sigma_{ess}(\wh A_0^+) = \max \sigma_{ess}(\wh A_0^-) =0
$$
and Theorem \ref{genau Def} concludes the proof.
\end{proof}

\section{Some classes of Sturm-Liouville operators} 
\label{s SomeCl}

As an illustration of the results from the previous sections,
we discuss in this section various potentials $q \in L^1_{loc} (\R)$ such that
%the differential expression \eqref{e ell}  is in limit point %case  at $+ \infty$ and at $-\infty$ and 
the
differential operator $A = (\sgn x) (-d^2/dx^2+q)$ is definitizable over specific subsets of $\overline{\C}$.
As before it is supposed that the differential expression \eqref{e ell}  is in limit point case  at $+ \infty$ and at $-\infty$ (for instance, the letter holds if \ $\liminf_{|x| \to \infty} \frac {q(x)}{x^2} > -\infty$, see e.g., \cite[Example 7.4.1]{Z05}).

\subsection{The case $q(x) \to - \infty$}

In this subsection we assume that for some $X>0$ the potential $q$ has the following properties on the interval $(X,+\infty)$: \begin{align} 
& q', \ q'' \ \text{exist and are continuous on} \ (X,+\infty ), \quad \quad q(x) < 0 , \quad  \quad q' (x) < 0 , \label{e qC2} \\ 
& %\text{and} \quad 
q''(x) \quad \text{is of fixed sign, \ 
i.e.,} \quad q'' (x_1) q'' (x_2) \geq 0 \quad \text{for all} 
\quad x_1,x_2 > X,  \label{e q''}\\ 
& \lim_{x \to +\infty} q(x) = -\infty , \quad 
\int_{X}^{+\infty} |q(x)|^{-1/2} dx = \infty, \quad 
\text{and}  \quad 
\limsup_{x \to +\infty} \frac{|q' (x)|}{|q(x)|^{p}} < \infty,  
\label{e q+inf}
\end{align}
where
$p \in (0,3/2)$ is a constant.
%\begin{gather*}
%q' (x) = O (|q(x)|^{p})  \qquad \text{as} \quad x \to +\infty.
%\end{gather*}
%\end{assumption}

Then the well-known result of Titchmarsh
(see e.g. \cite[Theorems 3.4.1 and 3.4.2]{LevSar90}) states that 
%under Assumption \ref{a q-inf} 
\eqref{e ell} 
%the differential expression $-d^2/dx^2+q$ 
is in the limit point case at
$+\infty$ and $\sigma (A_0^+) = \R$. Hence 
the set $\slr_A$ defined by \eqref{e slr} coincides with $\sigma_{ess} (A_0^-) \cup \infty$. By Theorem \ref{t DomLD},
there are two cases: 

\begin{description}
\item[(i)] Let $\sigma_{ess} (A_0^-) \neq \R$. Then the greatest domain over which $A$ is definitizable is $\Omega_A := \C \setminus \sigma_{ess} (A_0^-)$ (note that $\infty \not \in \Omega_A$). 
\item[(ii)] Let $\sigma_{ess} (A_0^-) = \R$. Then $\Omega_A  \cap \overline{\R} = \emptyset$ and there exists no domain $\Omega$ in $\overline{\C}$ such that $A$ is definitizable
over $\Omega$. In particular, the letter holds if the analogues of assumptions \eqref{e qC2}-\eqref{e q+inf} are fulfilled for $x \in (-\infty,0]$. 
\end{description}

\begin{example} \label{e x}
Let us consider the operator $A= (\sgn x )(-d^2/dx^2 - x)$.
By \cite[Theorem 6.6]{W87}
the differential expression $-d^2/dx^2 - x$ is in limit point case at $+\infty$ and $-\infty$.
Assumptions \eqref{e qC2}-\eqref{e q+inf} hold for $x \in (0,+\infty)$, hence $\sigma_{ess} (A_0^+) = \sigma (A) = \R$.
On the other hand, $\sigma_{ess} (A_0^-) = \emptyset$
(see Subsection \ref{ss q+inf} and \cite[Section 3.1]{LevSar90}). Therefore 
%Applying Theorems \ref{t q-inf} and \ref{t q+inf},
the operator $A$ is definitizable over $\C$ and
there exists no domain $\Omega$ in $\overline{\C}$ with $\infty \in \Omega$
such that $A$ is definitizable
over $\Omega$. By Proposition \ref{p nonreal}, 
the only possible accumulation point for non-real spectrum of $A$ is the
point $\infty$.
\end{example}

\subsection{The case $q(x) \to + \infty$ \label{ss q+inf}}

Let us assume that the following conditions 
holds with certain constants $X,c >0$:
%The potential $q$ is continuous on 
%semi-bounded from below on $\R_+$ 
\begin{gather} \label{e Molch}
q(x) \geq c \quad \text{for} \quad x>X, \quad 
\text{and for any} \ \omega>0, \quad 
\lim_{x \to +\infty} \int_{x}^{x+\omega} q(t) dt = +\infty .
\end{gather}
%for any $\omega>0$.
%Moreover, there exists $C,M >0$ such that
%\begin{equation} \label{vonunten}
%q(x) \geq -Cx^2 \quad \mbox{for } x < -M.
%\end{equation}
%\end{assumption}

%The semi-boundedness from below and (\ref{vonunten}) implies %that the differential expression
%$-d^2/dx^2 +q$ from Subsection \ref{definitionen} is in the %limit point case
%at $+\infty$ and at $-\infty$, see e.g.\ \cite[Theorem %6.6]{W87}.

%A.M. 
Mol\v{c}anov proved (see e.g., \cite[Lemma 3.1.2]{LevSar90} and \cite[Subsection 24.5]{Nai68II})
that %when Assumption \ref{a Molch} is fulfilled,
\eqref{e Molch} yields 
$\sigma_{ess} (A_0^+) = \emptyset$, i.e.,
the spectrum of the operator $A_0^+$ is discrete.
Besides, \eqref{e Molch} implies that $A_0^+$ is semi-bounded from below. 
It follows from the results of Subsection \ref{ss ldinf} that
\emph{the operator $A$ is definitizable over $\C$}.
More precisely,

\begin{description}
\item[(i)] Let the operator $A_0^-$ be semi-bounded from above. Then the operator
$A$ is definitizable, $\infty$ is a regular critical point of $A$ (cf. \cite{CL89}), and $A$ admits decomposition \eqref{e Adec}.

\item[(ii)] Let $A_0^-$ be not semi-bounded from above. Then
$A$ is definitizable over $\C$ and there exists no domain
$\Omega$ in $\overline{\C}$ with $\infty \in \Omega$ such that $A$ is definitizable over $\Omega$. The only possible accumulation point for non-real spectrum of $A$ is the
point $\infty$.
\end{description}

Note that $A_0^-$ is not semi-bounded from above if 
$\lim_{x \to -\infty} q(x) = - \infty$.

\subsection{Summable potentials}

%\begin{assumption} \label{a negsp}
%Assume that $A_0^+$ is semi-bounded from below, $A_0^-$ is semi-bounded from above,
%$\inf \sigma_{ess} (A_0^+) = 0$, and $\sup \sigma_{ess} (A_0^-) = 0$.
%\end{assumption}

We denote by $q_{neg}(x) :=\min \{q(x), 0 \}$, $x \in \R$.

\begin{assumption} \label{a qneg}
$ \displaystyle
\int_t^{t+1} |q_{neg} (x)| dx \to 0 $
%\qquad \text{
as
%} \quad
$|t| \to \infty.$
\end{assumption}
If Assumption \ref{a qneg} is fulfilled then
the differential expression $-d^2/dx^2 + q$ is in limit point case at $+\infty$ and $-\infty$, cf.\ \cite[Satz 14.21]{W03}. 
By \cite[Theorem 15.1]{W87}, $A_0^+$ is semi-bounded from below, $A_0^-$ is semi-bounded from above with
$$
\sigma_{ess} (A_0^+) \subset [0,+\infty) \quad
 \mbox{and} \quad \sigma_{ess} (A_0^-) \subset (-\infty,0].
$$
This implies that 
the negative spectrums of the operators $A_0^+$ and $-A_0^-$ consist of
eigenvalues,
\[
\sigma (\pm A_0^\pm) \cap (-\infty,0) = \{ \pm \la_n^\pm \}_1^{N^\pm} \subset \sigma_p (\pm A_0^\pm),
\]
where $0 \leq N^\pm \leq \infty$. Besides, $\lim_{n \to \infty} \la_n^\pm = 0$ if $N^\pm =\infty$.
Then, by Theorem \ref{genau Def}, $A$ is definitizable over
$\ov{\C} \setminus \{0\}$.
Theorems \ref{t DomLD}
and \ref{t RegInf} imply easily the following statement.
% Again, let $a_0^+$ and $a_0^-$
%be defined as in (\ref{a=a-}) and (\ref{a=a+}), respectively.

\begin{theorem} \label{p sum}
Let %$A = (\sgn x) (-d^2/dx^2 +q)$, and let
Assumption \ref{a qneg} be fulfilled.
Then the operator \\
$A = (\sgn x) (-d^2/dx^2 +q)$ admits the decomposition \eqref{e Adec}. Moreover, 
\begin{description}
\item[(i)] If  $\min \sigma_{ess} (A_0^+) >0$ or 
$\max \sigma_{ess} (A_0^-) <0$, then
$A$ is a definitizable operator and $\infty$ is a critical point of $A$.
\item[(ii)] If $\min \sigma_{ess} (A_0^+) = \max \sigma_{ess} (A_0^-) = 0$ and $N^+ + N^- < \infty $, then
$A$ is a definitizable operator, $0$ and $\infty$ are critical points of $A$ .
\item[(iii)] If $\min \sigma_{ess} (A_0^+) = \max \sigma_{ess} (A_0^-) = 0$ and $N^+ + N^- = \infty $,
then the operator $A$ is not definitizable. It is definitizable over
$\overline{\C} \setminus \{0\}$.
In particular, $0$ is the only possible
accumulation point of the non-real spectrum of $A$.
\end{description}
\end{theorem}

We mention (cf. \cite{B07}) that Assumption \ref{a qneg}, and therefore
the statements of Theorem \ref{p sum}, hold true if $q \in L^1 (\R)$.

\begin{remark}
By Theorem \ref{t RegInf} (see also \cite{CL89}) we have that  if the operator $A = (\sgn x) (-d^2/dx^2 +q)$ is definitizable, 
then $\infty$ is its regular critical point. 
In the case when $A$ has a finite critical point, the question of the character of this critical point is difficult (see \cite{CN96,CN98,FSh00,FadSh02,KarM04,KarM06,KarKos06} and references therein). Let us mention one case.
%
%\textbf{(i)}
Assume that $q$ is continuous in $\R$
and $\int_\R (1+x^2) |q(x)| dx <\infty$, then $\min \sigma_{ess} (A_0^+) = \max \sigma_{ess} (A_0^-) = 0$ and $N^+ < \infty$ and
$N^- < \infty$ (see \cite{LevSar90}
%\cite[Lemma 3.2.2]{LevSar90}
).
Therefore Theorem \ref{p sum} (as well as
\cite[Proposition 1.1]{CL89}) implies that $A= (\sgn x) (-d^2/dx^2+q)$ is definitizable. 
It was shown (implicitly) in \cite{FSh00} that $0$ is a regular critical point of $A$.
%
%\textbf{(ii)}In the case $q=6\frac{x^4-6|x|}{(|x|^3+3)^2}$ %(see \cite{KarKos06}), the point $0$ is a singular critical %point of the operator $A$ .
\end{remark}

%If $\lim_{x \to \infty} q(|x|) = 0$, then $\min \sigma_{ess} %(A_0^+) = \max \sigma_{ess} (A_0^-) = 0$ (see \cite{Z05}) 
%and therefore either statement (ii) or statement (iii) of %Theorem \ref{p sum} holds.
In the following case, more detailed information may be obtained.

\begin{corollary} \label{c sum2}
Suppose $\lim_{x \to \infty} q(|x|) = 0$. Then $\min \sigma_{ess} (A_0^+) = \max \sigma_{ess} (A_0^-) = 0$ 
and either the case (ii) or the case (iii) of Theorem \ref{p sum} takes place. Moreover, the following holds.
\begin{description}
\item[(i)] If \ $\liminf_{x \to \infty} x^2 q(|x|) > -1/4$, then
$A$ is a definitizable operator and $0$ and $\infty$ are
critical points of $A$.
\item[(ii)] If \ $\limsup_{x \to \infty} x^2 q(|x|) < -1/4$, then the operator $A$ is not definitizable. It is definitizable over $\overline{\C} \setminus \{0\}$.
\end{description}
\end{corollary}

\begin{proof}
The statement follows directly from \cite[Corollary XIII.7.57]{DSchV263}, which was proved in \cite{DSchV263} 
for infinitely differentiable $q$. Actually, this proof 
is valid for bounded potentials $q$. Finally, note that 
$\lim_{x \to \infty} q(|x|) = 0$ implies that $q$ is bounded 
on $(-\infty,-X]\cup [X,+\infty)$ with $X$ large enough.
On the other hand, $L^1$ perturbations of potential $q$ on any finite interval does not change $\sigma_{ess} (A_0^+)$, $\sigma_{ess} (A_0^-)$. Also such perturbations increase or decrease $N^+$, $N^-$ on finite numbers only due to Sturm Comparison Theorem (see e.g., \cite[Theorem 2.6.3]{Z05}). This completes the proof.
\end{proof}

\begin{example} \label{ex 2}
Let $q(x) = - \frac {1}{1+|x|}$. Then Corollary \ref{c sum2} yields that
the operator $A = (\sgn x) (-d^2/dx^2 +q)$ is not definitizable.
It is definitizable over $\overline{\dC} \setminus \{ 0\}$.
\end{example}

It was shown above that under certain assumption on the potential $q$ the operator \\
$A = (\sgn x) (-d^2/dx^2 +q)$ is not definitizable, but it is definitizable over the domain $\overline{\dC} \setminus \{ \lambda_0\}$, where $\la_0 \in \overline{\R}$ ($\lambda_0 = \infty$ in Example \ref{e x}  
and $\lambda_0 = 0$ in Example \ref{ex 2}). In this case, unusual spectral behavior may appear near points of the set $c(A)\cup \{ \la_0\}$ only ($c(A)$ is the set of critical points, see Subsection \ref{sectionDEFF}).
Indeed, a bounded spectral projection $E^A (\Delta)$ exists for any connected set $\Delta \subset \overline{\R} \setminus \{\la_0\}$ such that the endpoints of $\Delta$ do not belong to $c(A)\cup \{ \la_0\}$. Note also that $c(A)$ is at most countable and that 
$\la_0$ is the only possible accumulation point of the non-real spectrum of $A$. 

% Thus the following definition is natural. Assume that 
% \begin{description}
% \item[(i)] a self-adjoint (in Krein space sense) operator $T$ 
% is definitizable over the domain 
% $\Omega_\la \setminus \{\lambda\}$, where $\la \in \ov{\R}$ 
% and $\Omega_\la$ is an open neighbourhood of $\la$,
% \item[(ii)] $T$ is not definitizable over $\Omega_\la$.
% \end{description}
% Then $\la$ is said to be a \emph{quasi-critical point} of the 
% operator $A$.
% 
% The concept of regular critical point gives an approach to 
% the similarity problem due to the following fact:
% a nonnegative operator $T$ in a Krein space 
% $(\mathcal H, [\cdot,\cdot])$ is similar
% to a self-adjoint operator in the Hilbert space $\mathcal H$
% if and only if $0$ and $\infty$ are not singular critical 
% points of $T$ and $\ker (T) = \ker (T^2)$.
% (see e.g. \cite{CN96}).
% To  extend this method for nondefinitizable operators 
% considered in Section \ref{s SomeCl},
% one should carry the concept of regularity over super-critical points.

\section*{Acknowledgement}

The first author acknowledges the  hospitality and support of
the Technische Universit\"at Berlin and of the University of Zurich.

\vspace{0.5cm}

\noindent
Illya Karabash\\ 
Department of Partial Differential Equations\\
Institute of Applied Mathematics and Mechanics of NAS of Ukraine, R. Luxemburg str. 74\\
Donetsk 83114\\
Ukraine, \\
e-mail:
{\tt\small  karabashi@yahoo.com, karabashi@mail.ru}\\

\vspace{0.5cm}

\noindent
Carsten Trunk\\
Institut f\"ur Mathematik\\
Technische Universit\"at Berlin\\
Sekretariat MA 6-3\\
 Stra\ss e des 17.~Juni 136\\
 D-10623 Berlin\\
Germany, \\
e-mail: {\tt\small trunk@math.tu-berlin.de}

\end{document}